

\input amstex
\documentstyle{amsppt}
\nologo


\input xy
\xyoption{all}


\voffset=-5mm
\hsize=14.5cm
\vsize=22.5cm
\magnification=1100


\def\esp#1{\vskip#1cm}
\def\espaceavantpar{\esp{0.9}\vglue 1mm}
\def\espaceaprespar{\esp{0.2}}
\def\espaceavanttruc{\esp{0.4}}

\def\soupir{\kern2mm}


\newcount\parno         
\parno=0

\newcount\trucno        
\trucno=0

\newcount\bibleno       
\bibleno=0


\def\paragraph#1{\goodbreak
\espaceavantpar
\global\advance\parno by 1\trucno=0 \noindent{\bf \the\parno.\soupir{#1}}
\espaceaprespar}

\def\point{\global\advance\trucno by 1
\the\parno.\the\trucno.\soupir}

\def\truc#1{\espaceavanttruc\goodbreak
\noindent{\it\point #1.\rm} \ignorespaces}

\def\enonce#1{\goodbreak\espaceavanttruc\goodbreak
\noindent{\bf\point{#1}.}\soupir\ignorespaces\it}

\def\enonceplus#1#2{\goodbreak\espaceavanttruc\goodbreak
\noindent{\bf\point{#1}\soupir({#2}).}\soupir\ignorespaces\it}

\def\defi{\truc{Definition}}

\def\rem{\truc{Remark}}

\def\nota{\truc{Notation}}

\def\exo{\truc{Exercise}}

\def\exa{\truc{Example}}

\def\thm{\enonce{Theorem}}
\def\endthm{\rm}

\def\cor{\enonce{Corollary}}
\def\endcor{\rm}

\def\prop{\enonce{Proposition}}
\def\endprop{\rm}

\def\lem{\enonce{Lemma}}
\def\endlem{\rm}

\def\proof{\esp{0.2}\rm\noindent{\it Proof.\rm}\soupir\ignorespaces}
\def\qed{\hfill$\square$\goodbreak}


\def\bible#1;#2;#3;
{\medbreak \hangafter=1\hangindent=6mm \global\advance\bibleno by 1
\noindent\rlap{\the\bibleno.}\kern6mm{\ignorespaces #1},\ {\it\ignorespaces
#2},\ \ignorespaces #3.}


\def\disp{\displaystyle}
\def\isoto{\buildrel \sim\over\to}
\def\too{\mathop{\longrightarrow}}

\def\setminus{-}
\def\inv{^{-1}}
\def\sep{\,\big|\,}
\def\Lotimes{\otimes^L}
\def\unit{{1\!\!1}}

\def\suchthat{\text{\,s.t.\,}}
\def\red{{\text{red}}}
\def\ptred{{\text{pt.red}}}

\def\etal{{\sl et al.}}
\def\ie{{\sl i.e.}\ }

\def\N{\Bbb N}

\def\calo{{\Cal O}}

\def\Id{{\operatorname{Id}}}
\def\End{{\operatorname{End}}}

\def\PNil{\operatorname{PNil}}
\def\Kb{{\operatorname{K}}^{\operatorname{b}}}
\def\Db{{\operatorname{D}}^{\operatorname{b}}}
\def\Dperf{{\operatorname{D}}^{\operatorname{perf}}}
\def\Space{{\operatorname{Space}}}
\def\Spacered{{\operatorname{Space}_\ptred}}
\def\Spc{{\operatorname{Spc}}}
\def\Spec{{\operatorname{Spec}}}
\def\Supph{{\operatorname{Supph}}}


\def\satur#1{\left<{#1}\right>}
\def\bigsatur#1{\big<{#1}\big>}
\def\Bigsatur#1{\Big<{#1}\Big>}

\def\adh#1{\overline{\{#1\}}}
\def\adhs#1{\satur{#1}}

\def\ring#1{\calo_{#1}}
\def\ringred#1{\calo_{#1}'}


\def\xytriangle#1;#2;#3;#4;#5;#6;{
\xymatrix{#1\ar[r]^-{\disp #4}&#2\ar[r]^-{\disp #5}&#3\ar[r]^-{\disp
#6}&T(#1)}}

\def\magnxytriangle#1;#2;#3;#4;#5;#6;#7;{
\xymatrix@C=#1pt{#2\ar[rr]^-{\disp #5}&&#3\ar[rr]^-{\disp
#6}&&#4\ar[rr]^-{\disp #7}&&T(#2)}}



\def\namedef#1{\expandafter\def\csname #1\endcsname}
\def\newlabel#1#2{\global\namedef{r#1}{#2}}


\newlabel{defperf}{{2.1}{0}}
\newlabel{defthick}{{2.2}{0}}
\newlabel{deftenscat}{{2.3}{0}}
\newlabel{deftensthick}{{2.6}{0}}
\newlabel{thomclass}{{2.8}{0}}
\newlabel{hypnoeth}{{2.9}{0}}
\newlabel{thomplus}{{2.10}{0}}
\newlabel{consloc}{{2.11}{0}}
\newlabel{consIC}{{2.12}{0}}
\newlabel{thtrplus}{{2.13}{0}}
\newlabel{defatomic}{{3.3}{0}}
\newlabel{defatomicbis}{{3.4}{0}}
\newlabel{notaUF}{{3.6}{0}}
\newlabel{topology}{{3.8}{0}}
\newlabel{exunderly}{{3.9}{0}}
\newlabel{trispeclo}{{3.10}{0}}
\newlabel{defgeom}{{4.1}{0}}
\newlabel{defdense}{{4.2}{0}}
\newlabel{phibarbar}{{4.5}{0}}
\newlabel{propdefPhi}{{4.6}{0}}
\newlabel{suppressdense}{{4.7}{0}}
\newlabel{funcspc}{{4.11}{0}}
\newlabel{defJL}{{5.1}{0}}
\newlabel{presheaftenscats}{{5.3}{0}}
\newlabel{deftripresheaf}{{5.6}{0}}
\newlabel{defmolecular}{{5.8}{0}}
\newlabel{functrispc}{{5.11}{0}}
\newlabel{defEnd}{{6.1}{0}}
\newlabel{defpointnil}{{6.2}{0}}
\newlabel{kernil}{{6.3}{0}}
\newlabel{endequiv}{{6.4}{0}}
\newlabel{defringsp}{{6.6}{0}}
\newlabel{compare}{{6.9}{0}}
\newlabel{unionbar}{{7.1}{0}}
\newlabel{irredatomic}{{7.2}{0}}
\newlabel{underly}{{7.3}{0}}
\newlabel{Lfgeodense}{{7.4}{0}}
\newlabel{lemfinv}{{7.5}{0}}
\newlabel{undermorph}{{7.7}{0}}
\newlabel{thmOne}{{7.8}{0}}
\newlabel{folklore}{{8.1}{0}}
\newlabel{reconsSpec}{{8.3}{0}}
\newlabel{Dperfidele}{{8.4}{0}}
\newlabel{thmTwo}{{8.5}{0}}
\newlabel{maincorOne}{{8.6}{0}}
\newlabel{corGth}{{8.7}{0}}
\newlabel{trivexam}{{9.2}{0}}
\newlabel{Kblabla}{{9.3}{0}}
\newlabel{thmThree}{{9.7}{0}}
\newlabel{basc}{{9.8}{1}}
\newlabel{sga}{{9.8}{2}}
\newlabel{borl}{{9.8}{3}}
\newlabel{ega}{{9.8}{4}}
\newlabel{muka}{{9.8}{5}}
\newlabel{neem}{{9.8}{6}}
\newlabel{neek}{{9.8}{7}}
\newlabel{orlo}{{9.8}{8}}
\newlabel{poli}{{9.8}{9}}
\newlabel{rick}{{9.8}{10}}
\newlabel{thom}{{9.8}{11}}
\newlabel{thtr}{{9.8}{12}}
\newlabel{vakn}{{9.8}{13}}
\newlabel{verd}{{9.8}{14}}
\newlabel{weib}{{9.8}{15}}


\long\def\firstoftwo#1#2{#1}
\long\def\secondoftwo#1#2{#2}
\def\permute#1#2{\expandafter#2#1\null}

\def\tref#1{\expandafter\permute\csname r#1\endcsname\firstoftwo{}}
\def\bref#1{\expandafter\permute\csname r#1\endcsname\secondoftwo{}}
\def\cite#1{{\bf\bref{#1}}}
\def\citeplus#1{{\rm[{\bf\bref{#1}}]}}


\topmatter
\title
Presheaves of triangulated categories and reconstruction of schemes
\endtitle
\author
Paul Balmer
\endauthor
\affil SFB, University of M\"unster, Germany
\endaffil
\date
April 2002
\enddate
\address
\noindent Paul Balmer, SFB 478, Universit\"at M\"unster, Hittorfstr.\ 27, 48149
M\"unster, Germany
\endaddress
\email balmer\@math.uni-muenster.de
\endemail
\urladdr www.math.uni-muenster.de/u/balmer
\endurladdr
\abstract
To any triangulated category with tensor product $(K,\otimes)$, we associate a
topological space $\Spc(K,\otimes)$, by means of thick subcategories of $K$,
{\sl \`a la} Hopkins-Neeman-Thomason. Moreover, to each open subset $U$ of this
space $\Spc(K,\otimes)$, we associate a triangulated category ${\Cal K}(U)$,
producing what could be thought of as a {\it presheaf of triangulated
categories}. Applying this to the derived category
$(K,\otimes):=(\Dperf(X),\Lotimes)$ of perfect complexes on a noetherian scheme
$X$, the topological space $\Spc(K,\otimes)$ turns out to be the underlying
topological space of $X$; moreover, for each open $U\subset X$, the category
${\Cal K}(U)$ is naturally equivalent to $\Dperf(U)$.

As an application, we give a method to reconstruct any reduced noetherian
scheme $X$ from its derived category of perfect complexes $\Dperf(X)$,
considering the latter as a tensor triangulated category with $\Lotimes$.
\endabstract
\endtopmatter


\document
\rightheadtext{Triangulated categories and reconstruction of schemes}

\paragraph{Introduction}

Triangulated categories are ubiquitous. They appear in Homological Algebra,
Topology, Algebraic Geometry, Representation Theory, $K$-Theory, and so on. In
most cases, these triangulated categories come along with a tensor product, in
the very weak sense of Definition~\tref{deftenscat}, which is so flexible as to
admit the zero tensor product $\otimes=0$ in any case. Let us then rephrase our
anthem\,: {\it tensor} triangulated categories are ubiquitous.

In Algebraic Geometry, Grothendieck \etal~\citeplus{sga} defined the notion of
perfect complexes (see~\tref{defperf}) and introduced the triangulated category
$\Dperf(X)$, which is very close to the ``naive" derived category of bounded
complexes of vector bundles. Actually, $\Dperf(X)$ coincides with the latter on
a big class of schemes, including affine ones, but is in general more suitable
for geometry, as brightly demonstrated by Thomason in his state-of-the-art
treatment of algebraic $K$-theory of schemes~\citeplus{thtr}. Observe that this
triangulated category $\Dperf(X)$ is equipped with a tensor product
$\otimes=\Lotimes_{\calo_X}$.
\smallbreak

Our main result is\,:

\smallbreak
\noindent{\bf Theorem.} {\it The functor
$X\mapsto(\Dperf(X),\Lotimes_{\calo_X})$, from reduced noetherian schemes to
tensor triangulated categories, is faithful and reflects isomorphisms in the
following sense\,: given two reduced noetherian schemes, if there exists a
$\otimes$-equivalence of tensor triangulated categories between their derived
categories of perfect complexes then the two schemes are isomorphic.}
\smallbreak

This Theorem is to be found in Section~8. It follows from
Theorem~\tref{Dperfidele} and Corollary~\tref{maincorOne}. We also give there
partial results about non-necessarily reduced schemes.
\goodbreak

This result is interesting for various reasons. First of all, taken literally,
it says that the scheme invariant $X\mapsto\Dperf(X)$ taking values in tensor
triangulated categories is a ``complete invariant". It is a commonplace to say
that schemes are supposed to be more subtle mathematical objects than tensor
triangulated categories. Therefore, this result is rather surprising.

Secondly, we will actually define a left inverse functor to $X\mapsto\Dperf(X)$
on a big class of $\otimes$-triangulated categories. Therefore, almost
tautologically, some classical constructions which are performed on schemes,
can be {\it extended} to those tensor triangulated categories as well. We have
$K$-theory in mind. We do not claim that this is a clever way of defining
$K$-theory of triangulated categories, but we observe that the $K$-theory of a
scheme can be reconstructed from $\Dperf(X)$ and $\Lotimes$, since the whole
scheme can\,! This goes against the general feeling among $K$-theorists that
their invariant, or rather Quillen's, does not go through derived categories.
Further comment about this can be found below in Remark~\tref{Kblabla}.

Thirdly, the technique we develop in the proofs might perhaps be transposed to
other frameworks. Our backbone is one of the last articles of
Thomason's~\citeplus{thom}, where he classifies thick subcategories of
$\Dperf(X)$ in terms of $X$. Hopkins and Neeman (see~\citeplus{neem}) did the
case where $X$ is affine and noetherian without using the tensor product.
Thomason generalized it to non-affine, non-noetherian schemes, at the cost of
having to remember the tensor product $\Lotimes$ on $\Dperf(X)$. There is
already a considerable literature about ``classification of thick
subcategories" in more than one area of mathematics and we will not lose
ourselves in an overview of this material.

For our convenience, we introduce here the concept of {\it presheaf of
triangulated categories} (see~\tref{deftripresheaf}). This is again very
general and might be useful in other contexts as well. It allows us to speak of
things being true ``locally" in a triangulated category. Together with some
folklore observations about reconstructing reduced rings from their derived
category, this ``local" approach allows us to reconstruct reduced noetherian
schemes as well.

Our reconstruction holds more precisely for any scheme whose underlying space
is noetherian, see Remark~\tref{hypnoeth}. This is more general than being
noetherian, but the reader unfamiliar with this subtlety should just think of
noetherian schemes, like algebraic varieties for instance.

\smallbreak
Reconstruction results were already obtained by Bondal and
Orlov~\citeplus{borl}, without using the tensor structure, for smooth algebraic
varieties with ample either canonical or anticanonical sheaf. In that case,
$\Dperf(X)$ is equivalent to $\Db_{coh}(X)$, the derived category of coherent
$\calo_X$-modules. In general, it is impossible to reconstruct a variety from
$\Db_{coh}(X)$. There are known counterexamples to this due to Mukai, to Orlov
and to Polishchuk (see \citeplus{muka}, \citeplus{orlo} and \citeplus{poli}).

\medbreak
After Section~2 which contains some prerequisites, the article is divided into
two parts.

Part~I (Sections~3 to~6) gives the generalities about $\otimes$-triangulated
categories. We start with our definition of the {\it spectrum}
$\Spc(K,\otimes)$ of a $\otimes$-triangulated category. We insist on ``our"
since it is not clear that our approach would be exactly the right one in other
contexts, but it is good enough for our reconstruction purposes. In particular,
we give this spectrum a topology, in order to recover the underlying space of
$X$ when applied to $\Dperf(X)$. We also study the functoriality of this
construction $\Spc(-)$. We are led to consider what we call {\it geometric}
morphisms of $\otimes$-triangulated categories (Definition~\tref{defgeom}).
These include equivalences and morphisms of the form
$Lf^*:\Dperf(X)\to\Dperf(Y)$ for a morphism of noetherian schemes $f:Y\to X$.
Then we also consider presheaves on $\Spc(K,\otimes)$. The basic construction
being the {\it triangular presheaf} or presheaf of $\otimes$-triangulated
categories, denoted ${\Cal K}$, of Definition~\tref{defJL}. This uses the
general remark that the idempotent completion of a triangulated category can be
naturally triangulated, which is proved in~\citeplus{basc}. We then use rings
of endomorphisms locally to produce a presheaf of rings whose sheafification is
the wanted ringed space structure on $\Spc(K)$.

Part~II (Sections~7 to~9) is mostly devoted to the scheme case. We unfold and
apply the above definitions and constructions when $K$ equals $\Dperf(X)$. We
show that the underlying space of $X$ can be functorially recovered from
$(\Dperf(X),\Lotimes_X)$, which is tacitly in Thomason~\citeplus{thom}, and we
identify the presheaf of triangulated categories ${\Cal K}$ with
$U\mapsto\Dperf(U)$. We then check how the ringed space structure is linked to
the one of $X$. This leads us to the result. We end the paper with an example
and some further comments in Section~9.

\goodbreak
\medbreak
\noindent{\bf Acknowledgements:} The author sincerely thanks Pierre Deligne,
Bernhard Keller, Amnon Neeman, Christian Serp\'e and the referee for precious
comments or for interesting discussions.

\bigbreak
\goodbreak \noindent{\bf Table of Content}
\smallbreak
 1. Introduction.\par
 2. Preliminaries on triangulated categories.\par
\smallbreak\noindent Part I. Triangular presheaves\par
 3. Our definition of the spectrum.\par
 4. The geometric functors.\par
 5. A presheaf of tensor triangulated categories.\par
 6. Presheaves of rings.\par
\smallbreak\noindent Part II. Reconstruction of schemes.\par
 7. The underlying topological space.\par
 8. The main result.\par
 9. Comments.

\paragraph{Preliminaries on triangulated categories}

All categories under consideration are supposed to be essentially small.

An {\it abelian} or more generally an {\it exact} category is an additive
category equipped with a collection of so-called {\it exact sequences} which
encode the homological information. Similarly, a {\it triangulated} category is
an additive category $K$ where this homological information is encoded in {\it
exact triangles}. More precisely, $K$ comes equipped with a {\it translation},
\ie an endo-equivalence $T:K\to K$ and a class of {\it exact triangles} of the
form $\xytriangle P;Q;R;u;v;w;$, a.k.a.\ {\it distinguished triangles},
satisfying some axioms (see [\cite{verd}, II.1.1.1] or [\cite{weib},
Chapter~10]). A functor $\varphi:K\to K'$ is called {\it exact}, a.k.a.\ {\it
triangulated}, if it commutes with the translation up to isomorphism and
preserves exact triangles. Basic examples of triangulated categories are
derived categories of exact categories, stable categories of Frobenius
categories, Spanier-Whitehead categories of model categories, and so on. In
Algebraic Geometry and in Algebraic $K$-Theory, an important example is the
derived category $\Dperf(X)$ of perfect complexes on a scheme~$X$ (see
Grothendieck \etal~ [\cite{sga}, Expos\'e~I] or Thomason [\cite{thtr},
Section~2] and [\cite{thom}, 3.1]) that we briefly recall here.

\defi
Let $X$ be a scheme (below, $X$ will always be assumed to have a noetherian
underlying space). A {\it strict perfect complex} on $X$ is a bounded complex
of locally free ${\Cal O}_X$-modules of finite type. A complex $P$ of ${\Cal
O}_X$-modules is a {\it perfect complex} on $X$ if any point has a neighborhood
$U$ such that the restriction of $P$ to $U$ is quasi-isomorphic to some strict
perfect complex on $U$. We denote by $\Dperf(X)$ the full subcategory of
perfect complexes in $\text{D}(X)$, the derived category of the abelian
category of ${\Cal O}_X$-modules.

\espaceavanttruc
Thomason \citeplus{thom} gives a description of some subcategories of
$\Dperf(X)$ in terms of $X$ that we will use below. Let us recall what Verdier
\citeplus{verd} called an ``{\sl \'epaisse}" subcategory, and is here called a
``thick" subcategory. The following simplified version of the definition is due
to Rickard~\citeplus{rick}.

\defi
Let $K$ be a triangulated category. A subcategory $A\subset K$ is called {\it
thick} if it is a full triangulated subcategory (stable by isomorphisms,
translations, and taking cones) such that $P\oplus Q\in A$ with $P,Q\in K$
forces $P\in A$ and $Q\in A$.

\defi
We will call {\it tensor triangulated category} a triangulated category $K$
equipped with a covariant functor $\otimes:K\times K\to K$ which is exact in
each variable. A {\it morphism of $\otimes$-triangulated categories}, a.k.a.\
{\it$\otimes$-functor}, $\varphi:(K,\otimes)\to(L,\otimes)$ must commute with
$\otimes$ up to isomorphism.

\truc{Remarks}
\roster
\item
\hangafter=1\hangindent=5.5cm\rlap{$\smash{\xymatrix@C=21pt @R=25pt{T(P)\otimes
T(Q)\ar[d]_-{\disp\cong}\ar[r]^-{\disp\cong}&T(T(P)\otimes
Q)\ar[d]^-{\disp\cong}\\ T(P\otimes T(Q))\ar[r]^-{\disp\cong} &T^2(P\otimes
Q)}}$}\kern5.5cm Quite often, the diagram on the left is assumed to be
skew-commutative, as in the case of a derived tensor product. Here we will
assume commutativity {\it up to a sign}, which includes strict commutativity as
in Example~\tref{trivexam}.
\item
Observe that our definition of $\otimes$ is very flexible. We do not require
associativity or commutativity or anything of the like. Moreover, any
triangulated category can be equipped with a tensor structure\,: $\otimes=0$
for instance.
\endroster

\exa
Let $X$ be a scheme. The triangulated category $\Dperf(X)$ is equipped with
$\Lotimes$, the left derived functor of the usual tensor product
$-\otimes_{\Cal O_X}-$ (see if necessary [\cite{thom}, 3.1] and more references
there).

\defi
Let $(K,\otimes)$ be a $\otimes$-triangulated category. We will say that a
thick subcategory $A\subset K$ is {\it$\otimes$-thick} if $P\in A$ forces
$P\otimes Q\in A$ and $Q\otimes P\in A$ for all $Q\in K$.

\nota
Let $X$ be a scheme.
\roster
\item
Let $P\in\Dperf(X)$, denote by $\Supph(P):=\{x\in X\sep P_x\not\simeq 0\text{
in }\Dperf({\Cal O}_{X,x})\}$ the support of the homology of $P$.
\item
Let $Y\subset X$, denote by $\Dperf_Y(X):=\{P\in\Dperf(X)\sep\Supph(P)\subset
Y\}$ the full subcategory of $\Dperf(X)$ whose objects are those acyclic
outside $Y$.
\endroster

\enonceplus{Corollary}{from Thomason's classification}%
Let $X$ be a scheme with noetherian underlying space. Let ${\Cal C}$ be the
collection of $\otimes$-thick subcategories of $\Dperf(X)$ and let ${\Cal S}$
be the collection of specialization closed subsets of $X$, \ie subsets
$Y\subset X$ such that $y\in Y$ implies $\adh y\subset Y$. Then there are
inverse isomorphisms\,:
$$  \eqalign{\varphi:{\Cal C}&\to{\Cal S}\cr
            A&\mapsto \bigcup_{P\in A}\Supph(P)}
    \kern2cm
    \eqalign{\psi:{\Cal S}&\to{\Cal C}\cr
        Y&\mapsto\Dperf_Y(X)\vphantom{\bigcup_{P\in A}}}$$
which preserve inclusions.
\endcor

\proof
Thomason gives in [\cite{thom}, Theorem~3.15] a more general statement for
quasi-compact and quasi-separated schemes, in which the specialization closed
subsets of $X$ are replaced by unions of closed subsets whose complement is
quasi-compact. Under the noetherian assumption, this boils down to the above.
Recall that a scheme with noetherian underlying space is quasi-separated, see
[\cite{ega}, Section~6.1].
\qed

\truc{Remark and Definition}
Let us briefly comment on the hypothesis made throughout the article that our
schemes have noetherian underlying spaces. First of all, we need Thomason's
classification and thus we need our schemes to be quasi-compact and
quasi-separated. But later on, we want to work with all their open subschemes
and we will also need them to be quasi-compact and quasi-separated (actually
for different reasons, also going back to Thomason, see Theorem~\tref{thtrplus}
below). So we should deal with schemes $X$ such that\,:\hfill\break
\centerline{{\it $X$ is quasi-compact and quasi-separated, and
    any open subscheme of $X$ also is}\ \ ($\star$).}
This property is equivalent to the underlying topological space of $X$ being
noetherian. We will sometimes use the abbreviation {\it topologically
noetherian} for this property. It is an easy exercise to check that our
condition ($\star$) forces $X$ to be topologically noetherian. Conversely, if
$X$ is topologically noetherian, any open subscheme of $X$ is quasi-compact and
is moreover quasi-separated by [\cite{ega}, Corollary~6.1.13, p.~296].

\lem
Let $X$ be a topologically noetherian scheme.
\roster
\item
Let $Z\subset X$ be a closed subset. Then there exists a perfect complex
$P\in\Dperf(X)$ such that $\Supph(P)=Z$.
\item
Let $f:Y\to X$ be a morphism of topologically noetherian schemes and let
$P\in\Dperf(X)$. Then $\Supph(Lf^*(P))=f\inv(\Supph(P))$.
\endroster
\endlem

\proof
These two facts come respectively from Lemma~3.4 and Lemma~3.3 part (b) of
\citeplus{thom}. For more on $Lf^*$ see [\cite{thom}, 3.1].
\qed

\truc{Localization}
Let $K$ be a triangulated category and $J\subset K$ be a thick subcategory
(Definition~\tref{defthick}). We can define the {\it localization} or the {\it
quotient} $K/J$ as the triangulated category being universal for the following
property\,: there exists a functor $q:K\to L$ such that $q(P)\simeq0$ for all
$P\in J$. This category is obtained from $K$ via calculus of fractions by
inverting those morphisms of $K$ whose cone is in $J$. See more on that in
[\cite{verd}, Section~II.2]. Below we shall use for $K/J$ this model. The
category $K/J$ has the same objects as those of $K$ and morphisms are classes
of fractions
$$\mathop{\longleftarrow}_s\longrightarrow$$
where the cone of $s$ belongs to $J$. Clearly, when $K$ is moreover a
$\otimes$-triangulated category (Definition~\tref{deftenscat}) and when $J$ is
a $\otimes$-thick subcategory (Definition~\tref{deftensthick}), the
localization inherits a tensor product as well.

\truc{Idempotent completion}
A triangulated category $K$ is in particular additive and we can consider its
idempotent completion $\widetilde K$ as an additive category. It turns out that
$\widetilde K$ has a unique and natural structure of triangulated category such
that $K\to \widetilde K$ is exact. More about that can be found in
\citeplus{basc}. Below we shall always consider $\widetilde K$ to be built as
follows\,: objects are pairs $(P,p)$ where $p=p^2:P\to P$ is an idempotent and
morphisms $f:(P,p)\to (Q,q)$ are morphisms in $K$ such that $fp=qf=f$. Clearly,
when $K$ is moreover a $\otimes$-triangulated category then so is $\widetilde
K$.

\espaceavanttruc
Localization and idempotent completion allow the following result.

\thm
Let $X$ be a scheme with noetherian underlying space. Let $U\subset X$ be an
open subscheme and let $Z=X\setminus U$ its closed complement. Consider the
$\otimes$-thick subcategory $\Dperf_Z(X)$ of $\Dperf(X)$. Consider the quotient
$\Dperf(X)/\Dperf_Z(X)$ and its idempotent completion. Then there is a natural
equivalence of $\otimes$-triangulated categories\,:
$$\widetilde{\Dperf(X)/\Dperf_Z(X)}\too\Dperf(U)$$
which is compatible with the localization functors coming from $\Dperf(X)$.
\endthm

\proof
This is basically in Thomason-Trobaugh~\citeplus{thtr} and we only make it
explicit here for the convenience of the reader. It is clear that the
restriction functor $\Dperf(X)\to \Dperf(U)$ induces a functor\,:
$\Dperf(X)/\Dperf_Z(X)\too\Dperf(U).$ This functor is faithful by [\cite{thtr},
Proposition~5.2.4\,(a)] and full by [\cite{thtr}, Proposition~5.2.3\,(a)].
Observe that $\Dperf(U)$ is idempotent complete (exercise, see
also~\citeplus{basc}). It follows immediately that the idempotent completion of
$\Dperf(X)/\Dperf_Z(X)$ is also a full subcategory of $\Dperf(U)$. Now it
suffices to see that any object of $\Dperf(U)$ is a direct summand of the
restriction of an object of $X$. Let $F\in\Dperf(U)$, then $[F\oplus T(F)]=0$
in $K_0(U)$. The result now follows from [\cite{thtr}, Proposition~5.2.2\,(a)].
\qed

\espaceavantpar
\centerline{\bf Part I. Triangular presheaves}

\paragraph{Our definition of the spectrum}

\nota
Let $(K,\otimes)$ be a $\otimes$-triangulated category and $D$ be any
collection of objects in $K$. We will denote by $\satur{\,D\,}$ the smallest
$\otimes$-thick subcategory which contains $D$, that is the intersection of all
$\otimes$-thick subcategories which contain $D$. (See
Definition~\tref{deftensthick} for $\otimes$-thick.) As usual, when $D=\{d\}$
is reduced to a singleton, we shall write $\satur{d}$ instead of
$\satur{\{d\}}$.

\exo Let $(K,\otimes)$ be a $\otimes$-triangulated category and $D$ be any
collection of objects in $K$. Denote by $\Theta(D)$ the following subset of
objects of $K$\,: those $a\in K$ such that there exist $b,c\in D$, $e,e',f\in
K$, an exact triangle $\magnxytriangle 6;b;c;e\oplus e';;;;$ and an isomorphism
$a\simeq e$, $a\simeq e\otimes f$ or $a\simeq f\otimes e$. Let
$\Theta^0(D)=D\cup\{0\}$ and for each $n\geq 1$, let
$\Theta^n(D)=\Theta(\Theta^{n-1}(D))$. Show that $\satur{D}$ is the full
subcategory of $K$ on the following objects: $\cup_{n\in\N}\Theta^n(D)$.

\defi
Let $(K,\otimes)$ be a $\otimes$-triangulated category. A $\otimes$-thick
subcategory $A\subset K$ will be called {\it atomic} if the following condition
holds: whenever we have a collection of objects $D\subset K$ such that
$A\subset\satur{D}$ then there exists an element $d\in D$ such that
$A\subset\adhs{d}$.

\truc{Remarks}
Let $(K,\otimes)$ be a $\otimes$-triangulated category.
\roster
\item
An atomic $\otimes$-thick subcategory $A$ of $K$ is in particular {\it
principal} in the sense that there is an element $a\in A$ such that
$A=\adhs{a}$. The converse is false.
\item
It is an easy exercise to check that the condition of
Definition~\tref{defatomic} is equivalent to\,: whenever we have a collection
of $\otimes$-thick subcategories $\{B_i\subset K\,|\,i\in I\}$ such that
$A\subset\satur{\cup_{i\in I}B_i}$ then there exists an index $i\in I$ for
which $A\subset B_i$.
\endroster

\defi
Let $(K,\otimes)$ be a $\otimes$-triangulated category. We denote by
$\Spc(K,\otimes)$ the set of all non-zero atomic subcategories of $K$. We will
sometimes abbreviate this set by $\Spc(K)$ when the tensor structure is
understood. We now give it a topology.

\nota
Let $(K,\otimes)$ be a $\otimes$-triangulated category. For any object $a\in
K$, let
$$U(a):=\{B\in\Spc(K)\,\big|\,B\not\subset \adhs a\,\}.$$
We shall also make use later of the complement
$F(a):=\{B\in\Spc(K)\,\big|\,B\subset \adhs a\,\}$.

\lem
With the above notations, for any $a,b\in K$ we have $U(a)\cap U(b)=U(a\oplus
b)$. Moreover, $U(0)=\Spc(K)$.
\endlem

\proof
Observe that $\adhs{a\oplus b}=\satur{\,\{a,b\}\,}$. Thus by definition, any
atomic subcategory $B\subset K$ is contained in $\adhs{a\oplus b}$ if and only
if it is contained in $\adhs a$ or in $\adhs b$.
\qed

\defi
Taking the collection $\{U(a)\sep a\in K\}$ as a basis of a topology, we regard
$\Spc(K,\otimes)$ as a topological space. We call it {\it the spectrum} of the
$\otimes$-triangulated category $(K,\otimes)$.

\exa
Let $X$ be a topologically noetherian scheme. Then $\Spc(\Dperf(X),\Lotimes)$
is homeomorphic to the underlying topological space of $X$. This is
Theorem~\tref{underly} below.

\lem
Let $(K,\otimes)$ be a $\otimes$-triangulated category. Let $F\subset\Spc(K)$
be a closed subset and $A\in F$. For any $A'\in\Spc(K)$ such that $A'\subset
A$, we have $A'\in F$.
\endlem

\proof
For an elementary closed subset of the form $F=F(a)$ with $a\in K$
(see~\tref{notaUF}), it is immediate. Now any closed subset of $\Spc(K)$ is an
intersection of such subsets.
\qed

\paragraph{The geometric functors}

\espaceavanttruc
In order to make the construction $\Spc(-)$ functorial, we have to restrict
ourselves to some morphisms, that we call {\it geometric morphisms} of
$\otimes$-triangulated categories.

\defi
Let $\varphi:(K,\otimes)\to(L,\otimes)$ be a morphism of $\otimes$-triangulated
categories. (We denote both tensor products by $\otimes$ since no confusion can
occur from that.) We say that $\varphi$ is {\it geometric} if the following
condition is satisfied: For any collection of $\otimes$-thick subcategories
$\{C_i\sep i\in I\}$ in $K$ one has
$$\bigsatur{\varphi\big(\bigcap_{i\in I}C_i\big)}=\bigcap_{i\in I}\satur{\varphi(C_i)}.$$
Observe that the inclusion $\subset$ is always true. The above condition
roughly says that $\satur{\varphi(-)}$ commutes with arbitrary intersections.
\goodbreak

\defi
Let $\varphi:(K,\otimes)\to(L,\otimes)$ be a morphism of $\otimes$-triangulated
categories. We say that $\varphi$ is {\it dense} if $\satur{\varphi(K)}=L$. See
Remark~\tref{suppressdense}.

\exa
Let $f:Y\to X$ be a morphism of topologically noetherian schemes. Then
$Lf^*:\Dperf(X)\to\Dperf(Y)$ is geometric and dense. This is
Proposition~\tref{Lfgeodense} below.

\exa
Let $\varphi:(K,\otimes)\to(L,\otimes)$ be a $\otimes$-equivalence of
$\otimes$-triangulated categories. Then $\varphi$ is geometric and dense.

\lem
Let $\varphi: K\to L$ be a morphism of $\otimes$-triangulated categories. Let
$D\subset K$ be a collection of objects. Then
$\satur{\varphi(D)}=\satur{\,\varphi(\satur{D})\,}.$
\endlem

\proof
The left-hand side is clearly a subcategory of the right-hand side. Conversely,
the subcategory $\varphi\inv(\satur{\varphi(D)})$ of $K$ is $\otimes$-thick and
contains $D$, so it contains $\satur{D}$.
\qed

\truc{Proposition and Definition}
\it
Let $\varphi: K\to L$ be a morphism of $\otimes$-triangulated categories.
Assume that $\varphi$ is geometric and dense. Let $C\in\Spc(L)$. Define
$$\Phi(C):=\bigcap_{\disp H\subset K\text{ $\otimes$-thick subcategory}
\atop\disp \suchthat\ C\subset\satur{\varphi(H)}}H.$$ Then $\Phi(C)$ is a
non-zero atomic $\otimes$-thick subcategory of $K$. Moreover, we have
$C\subset\satur{\,\varphi(\Phi(C))\,}$.
\endprop

\proof
From $\varphi$ being geometric, we commute $\satur{\varphi(-)}$ with the above
intersection and get immediately
$$\satur{\varphi(\Phi(C))}=\bigcap_{\disp H\subset K\text{\ $\otimes$-thick subcategory}\atop
    \disp \suchthat\ C\subset\satur{\varphi(H)}}\satur{\varphi(H)}\kern5mm\supset C$$
which is the ``moreover part". Thus it is clear that $C\not=0$ implies
$\Phi(C)\not=0$.

Let $D\subset K$ be a collection of objects of $K$ such that
$\Phi(C)\subset\satur{D}$. Applying $\satur{\varphi(-)}$ to this relation,
using the previous Lemma and the above observation, we have that\,:
$C\subset\satur{\varphi(D)}$. Since $C$ is atomic, there is an element $d\in D$
with $C\subset\adhs{\varphi(d)}=\satur{\varphi(\adhs{d})}$. By the very
definition of $\Phi(C)$, this forces $\Phi(C)\subset\adhs{d}$ (think of
$H=\adhs{d}$\,). This gives the result.
\qed

\rem
We have only used $\varphi$ dense to be sure that there exists at least one
$\otimes$-thick subcategory $H$ of $K$ such that $\satur{\varphi(H)}\supset C$.
Pedantically speaking, requiring $\varphi$ to be geometric implies already that
$\varphi$ is dense. To see this, apply the defining condition to the empty
family\,: $\{C_i\sep i\in I\}=\emptyset$ in Definition~\tref{defgeom}.
Actually, we could even weaken the condition $\satur{\varphi(K)}=L$ into\,:
$\satur{\varphi(K)}\supset C$ for all $C\in\Spc(L)$.

Hereafter, we shall drop the mention ``dense" for geometric morphisms.

\prop
Let $\varphi:K\to L$ be a geometric morphism of $\otimes$-triangulated
categories. Let $\Phi:\Spc(L)\to\Spc(K)$ be defined as in~\tref{propdefPhi}.
Then $\Phi$ is continuous.
\endprop

\proof
Let $a\in K$. Recall the notation for closed subsets\,: $F(a)=\{B\in\Spc(K)\sep
B\subset\adhs a\,\}$ (see~\tref{notaUF}). It is enough to prove the following
Lemma\,:

\lem
With the above notations, $\Phi\inv(F(a))=F(\varphi(a))$.
\endlem

\proof
Let $C\in\Spc(L)$ be an atomic subcategory of $L$.

Assume that $C\in\Phi\inv(F(a))$. This means that $\Phi(C)\subset\adhs a$.
Applying $\satur{\varphi(-)}$ to both sides, using the ``moreover part" of
Proposition~\tref{propdefPhi} on the left, and Lemma~\tref{phibarbar} on the
right, we obtain $C\subset\adhs{\varphi(a)}$. This means $C\in F(\varphi(a))$.

Conversely assume that $C\in F(\varphi(a))$. This means
$C\subset\adhs{\varphi(a)}$. Using Lemma~\tref{phibarbar} again, this implies
that $H_0:=\adhs{a}$ is one of the $\otimes$-thick subcategories $H$ of $K$
satisfying $C\subset\satur{\varphi(H)}$. By the very definition of $\Phi(C)$,
we have $\Phi(C)\subset H_0=\adhs a$. This means $C\in\Phi\inv(F(a))$.

This finishes the proof of the Lemma and of the Proposition.
\qed

\defi
Let $\varphi:K\to L$ be a geometric morphism of $\otimes$-triangulated
categories. Define $\Spc(\varphi):\Spc(L)\to\Spc(K)$ to be the continuous map
$\Phi$ defined in~\tref{propdefPhi}.

\prop
We have a contravariant functor $\Spc(-)$ from $\otimes$-triangulated
categories with geometric morphisms to the category of topological spaces. Two
isomorphic geometric morphisms of $\otimes$-triangulated categories induce the
same map.
\endprop

\proof
Left as an easy familiarization exercise.
\qed

\exa
Let $f:Y\to X$ be a morphism of topologically noetherian schemes. Via the
identifications announced in Example~\tref{exunderly}, the morphism
$\Spc(Lf^*)$ coincides with the underlying morphism of spaces $f:Y\to X$. This
is Theorem~\tref{undermorph} below.

\paragraph{A presheaf of tensor triangulated categories}

\defi
Let $(K,\otimes)$ be a $\otimes$-triangulated category. Let
$V\subset\Spc(K,\otimes)$ be an open subset. Denote by
$J(V):=\satur{\,\bigcup_{A\in\Spc(K)- V}A\ }$. Since this is a $\otimes$-thick
subcategory, the idempotent completion of the localization
$${\Cal K}(V):=\widetilde{\big(K/J(V)\big)}$$
inherits a structure of $\otimes$-triangulated category (see~\tref{consloc}
and~\tref{consIC}). Below, we will denote by $q_V:K\to {\Cal K}(V)$ the natural
morphism.

\exa
Let $(K,\otimes)$ be a $\otimes$-triangulated category. Let $V=\Spc(K)$. Then
$J(V)=0$ and ${\Cal K}(V)=\widetilde{K}$ is the idempotent completion of $K$.

\prop
Let $(K,\otimes)$ be a $\otimes$-triangulated category. Let $V_1\subset V_2$ be
open subsets of $\Spc(K)$. Adopt the notations of~\tref{defJL}.
\roster
\item We have $J(V_2)\subset J(V_1)$.
\item There is a unique morphism of $\otimes$-triangulated categories $\rho_{V_1V_2}:{\Cal K}(V_2)\to {\Cal K}(V_1)$ such
that $q_{V_1}=\rho_{V_1V_2}\circ q_{V_2}$.
\endroster
Moreover, we have $\rho_{V_1V_3}=\rho_{V_1V_2}\circ\rho_{V_2V_3}$ for any
triple of open subsets $V_1\subset V_2\subset V_3$ in $\Spc(K)$.
\endprop

\proof
The first assertion is immediate. The second follows from the universal
property of localization and idempotent completion. Similarly for the
``moreover part". Actually, the above equalities are only isomorphisms of
functors. If we choose an explicit construction of all the localizations and
idempotent completions in terms of $K$, calculus of fractions and idempotents,
then these equalities hold strictly (see~\tref{consloc} and~\tref{consIC}).
\qed

\rem
Heuristically, we should consider $V\mapsto {\Cal K}(V)$ as a presheaf of
$\otimes$-triangulated categories on $\Spc(K,\otimes)$. We will not develop
here a whole theory of presheaves of triangulated categories, dealing with all
equalities which are in fact isomorphisms and the like. We leave these
improvements to the conscientious reader.

\exa
On a topologically noetherian scheme $X$, the above presheaf ${\Cal K}$
associated to the tensor triangulated category
$(K,\otimes)=(\Dperf(X),\Lotimes)$ produces $\Dperf(U)$ on each open $U\subset
X$. See Theorem~\tref{thmOne} below.

\defi
We shall call here {\it triangular presheaf}\,\ a pair $(X,{\Cal K})$
consisting of a topological space $X$ and a presheaf ${\Cal K}$ of triangulated
categories on $X$. Given a continuous map $f:X'\to X$ and a triangular presheaf
$(X',{\Cal K'})$, we denote as usual by $f_*{\Cal K'}$ the presheaf of
triangulated categories on $X$ defined by $f_*{\Cal K'}(V)={\Cal
K'}(f\inv(V)).$ A {\it morphism of triangular presheaves} $(f,{\Cal
F})\,:(X',{\Cal K'})\to(X,{\Cal K})$ consists of a continuous map $f:X'\to X$
and of a morphism of presheaves of triangulated categories ${\Cal F}:{\Cal
K}\to f_*{\Cal K'}$ on~$X$.

\rem
We now want to briefly study the functoriality of the above described
construction $K\mapsto (\Spc(K),{\Cal K})$ with respect to $K$, although,
strictly speaking, this is not necessary for the reconstruction of schemes. It
turns out that this functoriality is more tricky that one could expect and that
restricting ourselves to geometric morphisms is not enough. Therefore, we
restrict our study to a smaller class of triangulated categories, still
cautiously keeping in our basket all derived categories of perfect complexes on
topologically noetherian schemes.

\defi
Let $(K,\otimes)$ be a $\otimes$-triangulated category. We say that $K$ is {\it
molecular} if any $\otimes$-thick subcategory $B\subset K$ is generated by the
atomic subcategories it contains\,:
$$B=\Bigsatur{\bigcup_{\disp A\in\Spc(K)\atop\disp\suchthat A\subset B}A\ \ }.$$

\exa
Let $X$ be a topologically noetherian scheme. The $\otimes$-triangulated
category $\Dperf(X)$ is molecular. This is an easy exercise once we have
Theorem~\tref{underly} below.

\lem
Let $\varphi:K\to K'$ be a geometric morphism of $\otimes$-triangulated
categories. Assume that $K'$ is molecular. Let
$\Phi:=\Spc(\varphi):\Spc(K')\to\Spc(K)$ be the induced map
of~\tref{propdefPhi}. Let $V\subset\Spc(K)$. Recall the notations
of~\tref{defJL}. Then we have the following inclusion\,:
$$\varphi(J(V))\subset\ J\big(\Phi\inv(V)\big).$$
of subcategories of $K'$.
\endlem

\proof
Let $C\in\Spc(K')$ such that $C\subset \satur{\,\varphi(J(V))\,}$. We claim
that $\Phi(C)\in\Spc(K)\setminus V.$ It is obvious from the definition of
$\Phi$ given in~\tref{propdefPhi} that the inclusion
$C\subset\satur{\,\varphi(J(V))\,}$ implies $\Phi(C)\subset J(V)$. By
definition~\tref{defJL}, we have $J(V)=\satur{\cup_{A\in\Spc(K)\setminus V}A}$.
Now $\Phi(C)$ is atomic and is contained in $J(V)$, which forces the existence
of an $A\in\Spc(K)\setminus V$ with $\Phi(C)\subset A$. It follows from
Lemma~\tref{trispeclo} that $\Phi(C)$ also belongs to the closed subset
$\Spc(K)\setminus V$ of $\Spc(K)$. This proves the claim.

We have proved above that any $C\in\Spc(K')$ such that $C\subset
\satur{\,\varphi(J(V))\,}$ is contained in the $\otimes$-thick subcategory
$J\big(\Phi\inv(V)\big)$. We now use the hypothesis that $K'$ is molecular to
conclude that $\satur{\,\varphi(J(V))\,}$ is itself contained in
$J\big(\Phi\inv(V)\big)$.
\qed

\prop
Let $\varphi:K\to K'$ be a geometric morphism of $\otimes$-triangulated
categories. Assume that $K'$ is molecular. Let
$\Phi:=\Spc(\varphi):\Spc(K')\to\Spc(K)$ be the induced map and let ${\Cal K}$
and ${\Cal K'}$ be the respective presheaves of $\otimes$-triangulated
categories (see~\tref{defJL}). Then there is a morphism of presheaves of
$\otimes$-triangulated categories on $\Spc(K)$\,:
$${\Cal F}:{\Cal K}\too \Phi_*{\Cal K'}$$
whose global section ${\Cal F}(\Spc(K)):\widetilde K\to \widetilde K'$
coincides with the idempotent completion of~$\varphi$.
\endprop

\proof
Let $V\subset\Spc(K)$ be an open subset. Recall the notations of~\tref{defJL}.
From the above Lemma, we know that $\varphi$ induces a commutative (solid)
diagram\,:
$$\xymatrix{%
    0\ar[r]
    &J(V)\ar[r]\ar[d]
    &K\ar[r]\ar[d]_{\disp\varphi}
    &K/J(V)\ar[r]\ar@{-->}[d]
    &0
    \\
    0\ar[r]
    &J(\Phi\inv(V))\ar[r]
    &K'\ar[r]
    &K'/J(\Phi\inv(V))\ar[r]
    &0.
    }$$
Therefore $\varphi$ induces a morphism of $\otimes$-triangulated categories\,:
$K/J(V)\to K'/J(\Phi\inv(V))$. The idempotent completion of this morphism is
the wanted morphism ${\Cal K}(V)\to\Phi_*{\Cal K'}(V)$. The details are left to
the reader.
\qed

\cor
The construction $(K,\otimes)\mapsto (\Spc(K,\otimes),{\Cal K})$ described
in~\tref{defJL} defines a contravariant functor from molecular
$\otimes$-triangulated categories (see \tref{defmolecular}) and geometric
morphisms (see \tref{defgeom}) towards the category of triangular presheaves
(see \tref{deftripresheaf}).
\endcor

\proof
Obvious from the proof of the above Proposition.
\qed

\paragraph{Presheaves of rings}

\defi
Let $K$ be a triangulated category. We denote by $\End(K)$ the ring of those
endomorphisms of the identity functor on $K$ which commute with translation.
This is the same as the set of collections $(\alpha_a)_{a\in K}$ of
endomorphisms $\alpha_a:a\to a$ such that for any morphism $f:a\to b$ in $K$
one has $\alpha_b\,f=f\,\alpha_a$ and such that $\alpha_{T(a)}=T(\alpha_a)$.
Observe that $\End(K)$ is commutative.

\defi
An endomorphism of the identity $\alpha\in\End(K)$ is called {\it pointwise
nilpotent} if for any $a\in K$, $\alpha_a$ is nilpotent. We shall denote by
$\PNil(K)$ the ideal of pointwise nilpotent elements of $\End(K)$.

\prop
Let $K$ be a triangulated category and $D$ be a set of generators. (That is the
classical notion, \ie $K=\satur{D}$ for $\otimes=0$.) Assume that an element
$\alpha\in\End(K)$ is nilpotent on each object of $D$. Then $\alpha$ is
pointwise nilpotent.
\endprop

\proof
Clearly, if $\alpha_a^n=0$ then the same holds for a direct summand of $a$. So
it is enough to prove the following. Assume that
$$\xytriangle a;b;c;u;v;w;$$
is an exact triangle and that $\alpha_a$ and $\alpha_b$ are nilpotent, then so
is $\alpha_c$. Replacing $\alpha$ by some power of itself, we may as well
assume that $\alpha_a=0$ and $\alpha_b=0$. But then from
$\alpha_c\,v=v\,\alpha_b=0$ we deduce the existence of $\beta:T(a)\to c$ such
that $\alpha_c=\beta\,w$. It follows that
$\alpha_c^2=\beta\,w\,\alpha_c=\beta\,T(\alpha_a)\,w=0.$ More generally, in an
exact triangle of endomorphisms (obvious sense), if two of them are nilpotent
then so is the third.
\qed

\lem
Let $\varphi:K\simeq K'$ be an equivalence of triangulated categories. Then
$\varphi$ induces an isomorphism of commutative rings $\End(K)\simeq\End(K')$,
which preserves pointwise nilpotent elements.
\endlem

\proof
Let $\psi:K'\to K$ be an inverse equivalence and let
$\zeta:\Id_{K'}\isoto\varphi\psi$ be an isomorphism of functors. Define
$f:\End(K)\to\End(K')$ by the following formula\,:
$f(\alpha)_b=\zeta_b\inv\varphi(\alpha_{\psi(b)})\zeta_b$. Then the rest of the
proof is straightforward.
\qed

\lem
Let $K$ be a triangulated category and $J\subset K$ be a thick subcategory.
Consider the idempotent completion of the quotient $L=\widetilde{K/J}$ and
$q:K\to L$ the natural functor. Then there is a natural ring homomorphism
$\rho:\End(K)\to \End(L)$ such that $\rho(\alpha)_{q(a)}=q(\alpha_a)$ for any
$\alpha\in\End(K)$ and $a\in K$. Moreover, this homomorphism preserves the
pointwise nilpotent elements.
\endlem

\proof
Observe first of all that the previous Lemma allows us to choose our favorite
model for localization and idempotent completion, namely those
of~\tref{consloc} and~\tref{consIC}. The localization step is very easy. The
idempotent completion step is only easy\,: define $\rho(\alpha)$ on an object
$(a,p)$, where $p=p^2:a\to a$ to be
$\rho(\alpha)_{(a,p)}:=p\,\alpha_a=\alpha_a\,p$. The ``moreover part" is
obvious.
\qed

\defi
Let $(K,\otimes)$ be a $\otimes$-triangulated category. Consider the triangular
presheaf $(\Spc(K),{\Cal K})$ of Section~5. We consider a presheaf of rings
$V\mapsto\End({\Cal K}(V))$ on $\Spc(K,\otimes)$, whose sheafification will be
denoted by $\ring K$. We have defined a {\it ringed space}
$$\Space(K,\otimes):=(\Spc(K),\ring K).$$

Similarly, we consider the presheaf of rings $V\mapsto\End({\Cal
K}(V))/\PNil({\Cal K}(V))$ on $\Spc(K,\otimes)$, whose sheafification will be
denoted by $\ringred K$. We have defined a {\it ringed space}
$$\Spacered(K,\otimes):=(\Spc(K),\ringred K).$$

\exa
Let $X$ be a topologically noetherian scheme. We will see that the ringed space
$\Spacered\big(\Dperf(X)\big)$ is the reduced ringed space $X_\red$. In
general, $\Space\big(\Dperf(X)\big)$ will be a bigger ringed space than $X$.
See more in Theorem~\tref{thmTwo}.

\rem
It is unclear for what nice class of morphisms of $\otimes$-triangulated
categories the above constructions $\Space_{(\ptred)}(K,\otimes)$ could be made
functorial. At the very least, two equivalent categories have isomorphic ringed
spaces, by Lemma~\tref{endequiv}, which is enough for the reconstruction
purposes. For general morphisms $\varphi:K\to K'$, it seems rather
unpredictable when an endomorphism of the identity functor on $K$ extends to
$K'$.

\prop
Let $(K,\otimes)$ be a $\otimes$-triangulated category. There is a surjective
morphism $r_K:\ring K\too\ringred K$ of sheaves of rings on $\Spc(K,\otimes)$.
\endprop

\proof
The morphism is the sheafification of the obvious epimorphism
$\End\to\End/\PNil$.
\qed

\espaceavantpar
\centerline{\bf Part II. Reconstruction of schemes}

\paragraph{The underlying topological space}

\lem
Let $X$ be a topologically noetherian scheme. Let $\{Y_i\sep i\in I\}$ be a
collection of specialization closed subsets of $X$. Then their union
$Y=\cup_{i\in I}Y_i$ is also specialization closed and
$$\Dperf_Y(X)=\Bigsatur{\bigcup_{i\in I}\,\Dperf_{Y_i}(X)}.$$
\endlem

\proof
By Corollary~\tref{thomclass}, the $\otimes$-thick subcategory
$\bigsatur{\cup_{i\in I}\,\Dperf_{Y_i}(X)}$ of $\Dperf(X)$ has to be of the
form $\Dperf_Z(X)$ for some specialization closed subset $Z\subset X$. Since
Thomason's correspondence preserves inclusions, it is immediate that $Z$ must
be the smallest specialization closed subset containing all the $Y_i$'s, that
is their union.
\qed

\prop
Let $X$ be a topologically noetherian scheme and $Y\subset X$ be a
specialization closed subset. Then $\Dperf_Y(X)$ is a non-zero atomic
subcategory of $\Dperf(X)$ in the sense of Definition~\tref{defatomic} if and
only if $Y$ is non-empty, closed and irreducible.
\endprop

\proof
Assume that $Y$ is non-empty, closed and irreducible. There is a point $y\in Y$
such that $Y=\adh y$. Let $A:=\Dperf_Y(X)$ in $K:=\Dperf(X)$. First of all
$A\not=0$ by Lemma~\tref{thomplus}\,(1). We have to see that $A$ is atomic. We
choose to verify Definition~\tref{defatomic} via the equivalent condition of
Remark~\tref{defatomicbis}\,(2). Choose a collection of $\otimes$-thick
subcategories $\{B_i\subset K\,|\,i\in I\}$ such that
$A\subset\satur{\cup_{i\in I}B_i}$. For each $B_i$ take the unique
specialization closed subset $Y_i\subset X$ such that $B_i=\Dperf_{Y_i}(X)$.
The assumption that $A\subset\satur{\cup_{i\in I}B_i}$ and the previous Lemma
imply that $Y\subset\cup_{i\in I}Y_i$. Since $y\in Y$, this forces $y\in Y_i$
for some $i\in I$; therefore $Y=\adh y \subset Y_i$ and $A\subset B_i$.

Conversely assume that $A:=\Dperf_Y(X)$ is non-zero and atomic. From $A$
non-zero one has $Y$ non-empty. Since $Y\subset\cup_{y\in Y}\adh y$, the
previous Lemma forces $A\subset\bigsatur{\cup_{y\in Y}\Dperf_{\adh y}(X)}$.
Using that $A$ is atomic, there exists a $y\in Y$ such that
$A\subset\Dperf_{\adh y}(X)$. This in turn implies that $Y\subset\adh y$. That
is $Y=\adh y$ is an irreducible closed subset of $X$.
\qed

\thm
Let $X$ be a scheme with noetherian underlying space. Consider the map
$$\eqalign{E:\ X&\too \Spc\big(\Dperf(X)\big)\cr
x&\mapsto\Dperf_{\overline{\{x\}}}(X)=
    \big\{P\in\Dperf(X)\sep\Supph(P)\subset\overline{\{x\}}\,\big\}.}$$
This is a well-defined homeomorphism between the underlying space of $X$ and
$\Spc\big(\Dperf(X)\big)$ with the topology of Definition~\tref{topology}.
\endthm

\proof
The fact that $E:X\to \Spc(\Dperf X)$ is well defined and bijective follows
immediately from Proposition~\tref{irredatomic} and Thomason's classification
(Corollary~\tref{thomclass}). We want to prove that $E$ respects the topology.
Recall from Definition~\tref{topology}, that the topology on $\Spc(\Dperf X)$
is given by the following basis ${\Cal B}:=\{U(a)\sep a\in\Dperf (X)\}$, where
$U(a)=\{B\in\Spc(\Dperf X)\sep B\not\subset \adhs a\,\}$. Let $Y=\Supph(a)$; it
is a closed subset of $X$ and we have $\Dperf_{Y}(X)=\satur{a}$. It is clear
that $\Spc(K)- U(a)=\{E(x)\sep x\in X, \suchthat \Dperf_{\adh x}(X)\subset
\Dperf_Y(X)\}= \{E(x)\sep x\in X, \suchthat \adh x\subset Y
\}=E(Y)=E(\Supph(a))$.

In other words, under our bijection $E$, the basis ${\Cal B}$ of $\Spc(\Dperf
X)$ corresponds to the following collection of open subsets of $X$: ${\Cal
B}'=\{X- \Supph(a)\sep a\in\Dperf(X)\}.$ Conversely, since $X$ is topologically
noetherian, any open in $X$ is an element of ${\Cal B}'$, by
Lemma~\tref{thomplus}\,(1).
\qed

\prop
Let $f:Y\to X$ be a morphism of topologically noetherian schemes. Let
$\varphi:=Lf^*:\Dperf(X)\to\Dperf(Y)$. Then $\varphi$ is geometric and dense in
the sense of Definitions~\tref{defgeom} and ~\tref{defdense}.
\endprop

\lem
With the notations of the above Proposition, let $Z\subset X$ be specialization
closed. Observe that $f\inv(Z)$ is specialization closed in $Y$. We have
$$\Bigsatur{\varphi(\Dperf_Z(X))}=\Dperf_{f\inv(Z)}(Y).$$
\endlem

\proof
Let $y\in f\inv(Z)$. By continuity of $f$, we have $f(\adh y)\subset
\adh{f(y)}\subset Z$ since $Z$ is specialization closed. This means that
$\adh{y}\subset f\inv(Z)$. Thus $f\inv(Z)$ is specialization closed.

Therefore, the right-hand category is a $\otimes$-thick subcategory of
$\Dperf(Y)$. The following inclusion\,:
$\bigsatur{\varphi(\Dperf_Z(X))}\subset\Dperf_{f\inv(Z)}(Y)$ is a direct
consequence of Lemma~\tref{thomplus}\,(2). We want to prove the converse
inclusion. By Thomason's classification, there exists $W$ specialization closed
in $Y$ such that $\bigsatur{\varphi(\Dperf_Z(X))}=\Dperf_W(Y)$. We are left to
prove $W\supset f\inv(Z)$.

Let $x\in Z$. Then there is a complex $e_x\in\Dperf(X)$ such that
$\Supph(e_x)=\adh x$ by Lemma~\tref{thomplus}\,(1). Since
$\varphi(e_x)\in\bigsatur{\varphi(\Dperf_Z(X))}$ we have
$W\supset\Supph(\varphi(e_x))=f\inv(\Supph(e_x))$ by
Lemma~\tref{thomplus}\,(2). This means that $W\supset f\inv(\adh x)\supset
f\inv(\{x\})$ and this holds for arbitrary $x\in Z$. In other words, we have
$W\supset f\inv(Z)$ which is the claim.
\qed

\truc{Proof of Proposition~\tref{Lfgeodense}}
Let $\{C_i\sep i\in I\}$ be a collection of $\otimes$-thick subcategories in
$\Dperf(X)$ and, for each $i\in I$, let $Z_i\subset X$ be the specialization
closed subset such that $C_i=\Dperf_{Z_i}(X)$. Let also $Z:=\cap_{i\in I}Z_i$
and observe that $f\inv(Z)=\cap_{i\in I}f\inv(Z_i)$. We then compute directly,
using the above Lemma:
$$\eqalign{\Bigsatur{\varphi(\bigcap_{i\in I}C_i)}
    &=\Bigsatur{\varphi(\bigcap_{i\in I}\Dperf_{Z_i}(X))}
    =\bigsatur{\varphi(\Dperf_Z(X))}=\Dperf_{f\inv(Z)}(Y)=
    \cr
    &=\bigcap_{i\in I}\Dperf_{f\inv(Z_i)}(Y)
    =\bigcap_{i\in I}\bigsatur{\varphi(\Dperf_{Z_i}(X))}
    =\bigcap_{i\in I}\bigsatur{\varphi(C_i)}.}$$
This proves that $\varphi$ is geometric. It is dense because
$\calo_Y\simeq\varphi(\calo_X)\in\bigsatur{\varphi(\Dperf(X))}$.
\qed

\thm
Let $f:Y\to X$ be a morphism of schemes with noetherian underlying spaces.
Consider the morphism of $\otimes$-triangulated categories
$Lf^*:\Dperf(X)\to\Dperf(Y)$, that we know to be geometric. Then the following
diagram of topological spaces\,:
$$\xymatrix@C=60pt @R=30pt{
    Y\ar[r]^{\disp f}\ar[d]_{\disp E}^{\disp\simeq}&X\ar[d]_{\disp E}^{\disp\simeq}
    \\
    \Spc\big(\Dperf(Y)\big)\ar[r]_{\disp\Spc(Lf^*)}
    &\Spc\big(\Dperf(X)\big)
    }$$
is commutative. The vertical homeomorphisms are those of
Theorem~\tref{underly}, the lower line involves the functor $\Spc(-)$ of
Section~4.
\endthm

\proof
Let $y\in Y$. Let us use the following notations\,: $K:=\Dperf(X)$,
$L:=\Dperf(Y)$, $\varphi:=Lf^*: K\to L$ and $\Phi:=\Spc(\varphi):
\Spc(L)\to\Spc(K)$. We want to compute $\Phi(E(y))$. By
Definition~\tref{propdefPhi},
$$\Phi(E(y))=\bigcap_{\disp H\subset K\text{\ $\otimes$-thick subcategory}
\atop\disp \suchthat\ E(y)\subset\satur{\varphi(H)}}H.$$ Using Thomason's
classification to replace the $H$'s and using Lemma~\tref{lemfinv}, it is easy
to see that $\Phi(E(y))=\Dperf_W(X)$ where $W\subset X$ is given by
$$W=\bigcap_{\disp Z\subset X\text{\ specialization closed}\atop\disp \suchthat\ y\in f\inv(Z)}Z.$$
In other words, $W$ is the smallest specialization closed subset of $X$ which
contains $f(y)$. This is $\adh{f(y)}$. This means that
$\Phi(E(y))=\Dperf_{\adh{f(y)}}(X)=E(f(y))$, which is the claim.
\qed

\thm
Let $X$ be a scheme with noetherian underlying space. Let $U$ be any open
subscheme and $Z=X-U$ its closed complement. Let $V=E(U)$ be the corresponding
open in $\Spc\big(\Dperf(X)\big)$. Then, in the notations of
Definition~\tref{defJL}, we have $J(V)=\Dperf_Z(X)$ and there is a natural
equivalence of tensor triangulated categories with unit\,:
$${\Cal K}(V)\simeq\Dperf(U)$$
in a coherent way with respect to inclusions $U_1\subset U_2$ -- see
Proposition~\tref{presheaftenscats}.
\endthm

\proof
By Definition~\tref{defJL} and Theorem~\tref{underly},
$J(V)=\bigsatur{\cup_{x\in Z}\Dperf_{\adh x}(X)}$. From Lemma~\tref{unionbar},
we have then $J(V)=\Dperf_Z(X)$ since $\cup_{x\in Z}\adh x=Z$. The second
assertion follows immediately from Theorem~\tref{thtrplus}.
\qed

\paragraph{The main result}

\espaceavanttruc
\noindent We will start with a folklore observation. Let us first of all recall
that over a commutative ring $R$ the derived category $\Dperf(R)$ is simply the
homotopy category $\Kb({\Cal P}(R))$ of bounded complexes of finitely generated
projective $R$-modules.

\prop
Let $R$ be a commutative ring. Consider the ring homomorphism
$\lambda:R\to\End(\Dperf(R))$ sending $r\in R$ to the multiplication by $r$
everywhere (see~\tref{defEnd}).
\roster
\item
This homomorphism admits a left inverse $\sigma:\End(\Dperf(R))\to R$ given by
$\alpha\mapsto\alpha_\unit$, where $\unit=(\cdots\to0\to0\to
R\to0\to0\to\cdots)$ is concentrated in degree zero.
\item
The kernel of $\sigma$ is made of pointwise nilpotent elements
(see~\tref{defpointnil}).
\item
These homomorphisms induce isomorphisms
$R_\red\simeq\End(\Dperf(R))/\PNil(\Dperf(R))$.
\endroster
\endprop

\proof
The first assertion is obvious. The second follows from the fact that $\unit$
is a generator of $\Dperf(R)$ and from Proposition~\tref{kernil}. Now $\lambda$
and $\sigma$ clearly induce homomorphisms on the level of $R_\red$ and
$\End(\Dperf(R))/\PNil(\Dperf(R))$. The third assertion follows from (1) and
(2).
\qed

\rem
The above can be generalized to non-commutative rings, reconstructing in this
way the reduced ring of the center of $R$.

\cor
Let $R$ be a commutative ring. The topological space $\Spec(R)$ with its
Zariski topology is an invariant of the derived category $\Dperf(R)$, even for
$R$ non-noetherian.
\endcor

\proof
This is immediate from $\Spec(R_\red)=\Spec(R)$ and the above Proposition. The
construction $\End(-)/\PNil(-)$ was entirely made on the level of triangulated
categories.
\qed

\thm
Let $f,g: Y\to X$ be two morphisms of schemes and suppose that the two induced
morphisms of triangulated categories, $Lf^*$ and $Lg^*\,:\
\Dperf(X)\to\Dperf(Y)$ are isomorphic. Then $f=g$.
\endthm

\proof
We know already from Theorem~\tref{undermorph} that the underlying morphisms of
topological spaces $f$ and $g$ are equal (see also Proposition~\tref{funcspc}).
The end goes easily, for instance as follows.

Assume first that $X=\Spec(A)$ and $Y=\Spec(B)$ are affine and that $f,g\,:A\to
B$ are ring homomorphisms such that $\,Lf^*\simeq Lg^*$. Evaluating this
isomorphism $\,Lf^*\simeq Lg^*$ on the object $A$, viewed in $\Dperf(A)$ as a
complex concentrated in degree zero, we get an isomorphism of $B$-modules
$B\simeq B$. This is just multiplication by an invertible element $s\in B$.
Now, we have two commutative diagrams\,:
$$\xymatrix@C=40pt{A\ar[r]^-{\disp f}\ar[d]_-{\disp \mu}^-{\disp\simeq}
    &B\ar[d]^-{\disp \mu}_-{\disp\simeq}
    \\
    \End_{\Dperf(A)}(A)\ar[r]_-{\disp Lf^*}
    &\End_{\Dperf(B)}(B)
    }\kern5mm\text{and}\kern5mm
    \xymatrix@C=40pt{A\ar[r]^-{\disp g}\ar[d]_-{\disp \mu}^-{\disp\simeq}
    &B\ar[d]^-{\disp \mu}_-{\disp\simeq}
    \\
    \End_{\Dperf(A)}(A)\ar[r]_-{\disp Lg^*}
    &\End_{\Dperf(B)}(B)
    }$$
where the isomorphisms $\mu$ just send an element to the multiplication by it.
Our assumption that $Lf^*\simeq Lg^*$ boils down to $s\cdot f(a)\cdot
s\inv=g(a)$ for any $a\in A$. This means $f=g$.

Consider now $X$ and $Y$ not necessarily affine. Given $U\subset X$ and
$V\subset f\inv(U)=g\inv(U)$, with $U$ and $V$ both affine, we want to prove
that $Lf^*\simeq Lg^*$ implies $L\tilde f^*\simeq L\tilde g^*$, where $\tilde
f$ and $\tilde g$ are the restrictions of $f$ and $g$ to morphisms $V\to U$, as
in the left-hand diagram below\,:
$$\xymatrix{V\ar[r]\ar[rd]_{\disp\tilde f}
    &{f\inv}(U)\ar[r]\ar[d]
    &Y\ar[d]^{\disp f}
    \\
    &U\ar[r]
    &X
    }\kern2cm
    \xymatrix{\Dperf(X)\ar[d]_-{\disp q}\ar[r]^-{\disp Lf^*}
    &\Dperf(Y)\ar[d]^-{\disp q}
    \\
    \Dperf(U)\ar[r]_-{\disp L\tilde f^*}
    &\Dperf(V)
    }$$
We also have the above right-hand commutative diagram, where we denote by $q$
the restrictions and similarly for $g$ instead of $f$, {\sl mutatis mutandis}.
Now recall from Theorem~\tref{thtrplus}, that $\Dperf(U)$ is just the
idempotent completion of a localization of $\Dperf(X)$. Thus $L\tilde f^*$ is
characterized up to isomorphism by $L\tilde f^*\circ q$. Therefore, since
$Lf^*\simeq Lg^*$ forces $L\tilde f^*\circ q\simeq q\circ Lf^*\simeq q\circ
Lg^*\simeq L\tilde g^*\circ q$, we have $L\tilde f^*\simeq L\tilde g^*$. The
first part of the proof implies that $\tilde f=\tilde g$. This being true for
arbitrary affine open $U\subset X$ and $V\subset f\inv(U)=g\inv(U)$, it forces
$f=g$.
\qed

\thm
Let $X$ be a scheme with noetherian underlying space. Recall the notations of
Definition~\tref{defringsp}. There is a natural commutative diagram of sheaves
of rings on the space $X$\,:
$$\xymatrix{%
    \calo_X\ar[r]\ar[d]
    &\calo_{\Dperf(X)}\ar[d]^-{\disp r_{\Dperf(X)}}
    \\
    \calo_X^\red\ar[r]
    &\calo_{\Dperf(X)}'}$$
where the vertical homomorphisms are the natural epimorphisms
(see~\tref{compare} for the right one). Moreover, the upper morphism is split
injective and the lower morphism is an isomorphism. That is, there is a natural
ringed space isomorphism $X_\red\simeq\Spacered(\Dperf(X),\Lotimes)$.
\endthm

\proof
Observe first of all that those four sheaves of rings are defined on the same
topological space $\Spc(\Dperf(X),\Lotimes)=X$, via the identifications of
Theorem~\tref{underly}.

For each open $U\subset X$ there is a ring homomorphism
$\lambda:\calo_X(U)\to\End(\Dperf(U))$ given by the multiplication. This is
compatible with localization and idempotent completion and induces a morphism
of presheaves of rings $\lambda:\calo_X(-)\to\End({\Cal K}(-))$. Its
sheafification is the wanted morphism $\calo_X\to\ring{\Dperf(X)}$.

Conversely, there is a morphism $\sigma:\End(\Dperf(U))\to\calo_X(U)$ given by
evaluation at $\calo_U$ when $U$ is affine. This is again compatible with
restriction between affine open subsets. Since the latter form a sub-basis of
$X$, we have the wanted section $\sigma:\calo_{\Dperf(X)}\to\calo_X$.

Of course, nilpotent elements in $\calo_X(U)$ go to pointwise nilpotent
elements of $\End(\Dperf(U))$ and we have therefore the announced commutative
diagram. To see that the lower morphism $\lambda:\calo_{X_\red}(-)\to\End({\Cal
K}(-))/\PNil({\Cal K}(-))$ is an isomorphism, it suffices to do it locally.
This is Proposition~\tref{folklore}\,(3).
\qed

\cor
Let $X$ and $Y$ be topologically noetherian schemes. Assume that there is an
equivalence of $\otimes$-triangulated categories
$(\Dperf(X),\Lotimes)\simeq(\Dperf(Y),\Lotimes)$, not necessarily coming from a
morphism between $X$ and $Y$. Then the reduced schemes $X_\red$ and $Y_\red$
are isomorphic.
\endcor
\qed

\cor
Let $X$ and $Y$ be noetherian schemes. Assume that there is an equivalence of
$\otimes$-triangulated categories
$(\Dperf(X),\Lotimes)\simeq(\Dperf(Y),\Lotimes)$. Then $X$ and $Y$ have the
same $G$-theory ($K$-theory of coherent $\calo_X$-modules).
\endcor

\proof
We know by {\sl d\'evissage} that $G$-theory does not distinguish $X$ from
$X_\red$.
\qed

\paragraph{Comments}

\rem
We do not claim that our definition of $\Spc(-)$ will be useful as it stands
for triangulated categories not necessarily coming from schemes via
$\Dperf(-)$. The first dissuasion comes from non-noetherian schemes, even in
the affine case. See for instance Neeman's comments in [\cite{neem}, Section~4
and B\"okstedt's appendix]. Anyway, this question is interesting.

\exa
Consider one of the simplest examples of tensor triangulated category with
unit. Let $k$ be a field. Let $K=k-mod$ be the category of finite dimensional
vector spaces, let $T=\Id$ (!) and call a triangle $V_1\to V_2\to V_3\to V_1$
exact when the (long) sequence of vector spaces is exact at $V_1$, $V_2$ and
$V_3$. This is a triangulated category. It is equipped with tensor
$\otimes=\otimes_k$. Now let $A$ be a $\otimes$-thick subcategory of $K$.
Either $A=0$ or it contains some sum of copies of $k$ and therefore it contains
$k$ since it is thick, and thus $A=K$. In particular, $\Spc(K)=\{K\}$ is a
point. Similarly, $\Spc(\Dperf(k))=\Spec(k)$ is reduced to a point. In both
cases, the sheaf of rings is simply $k$. This shows that the above ringed space
constructions lose a lot of information about triangulated categories, which is
of course not a surprise.

\rem
Corollary~\tref{corGth} suggests that $K$-theory might be defined on the level
of tensor triangulated categories. Note that a good $K$-theory of triangulated
categories is still to be found. (If ``good" supposes agreement with Thomason's
$K$-theory on schemes and some long exact sequences on the triangular level.)
Of course, there is Thomason-Trobaugh~\citeplus{thtr}, but the second part of
the title ``Higher algebraic $K$-theory of schemes {\it and of derived
categories}" is slightly misleading. The $K$-theory ``of $\Dperf(X)$" requires
actually more data than just $\Dperf(X)$, namely it requires the underlying
bi-Waldhausen category of perfect complexes and quasi-isomorphisms. Similarly,
more recent attempts by Neeman~\citeplus{neek} to define $K$-theory of
triangulated categories still faces some drawbacks, including the recourse to
underlying models, leading in particular to his construction not being
functorial. Vaknin~\citeplus{vakn} has an improved functorial version, but
still uses models. The hope resulting from the above reconstruction result is
that models could be replaced, for $K$-theory purposes, by tensor products. In
the author's opinion this would be a conceptual and esthetic improvement.

\rem
There is another way of reconstructing the structure sheaf if we do not want to
assume that the scheme is reduced. Namely, we can keep track of the {\it unit}
for the tensor product $\unit=\calo_X$. In other words, we keep more
information and consider the functor $X\mapsto(\Dperf(X),\Lotimes,\calo_X)$ as
going from topologically noetherian schemes to tensor triangulated categories
with unit $(K,\otimes,\unit)$. For this, we take the following definition of a
unit.

\defi Given a $\otimes$-triangulated category $(K,\otimes)$, we will call an element $\unit\in K$
a {\it unit for $\otimes$} if there are natural isomorphisms
$\eta_P:P\otimes\unit\simeq P$ and $\mu_Q:\unit\otimes Q\simeq Q$ for all
$P,Q\in K$. Here, we shall also assume that
$\eta_\unit=\mu_\unit:\unit\otimes\unit\to\unit$. A morphism of tensor
triangulated categories with unit must respect the unit up to isomorphism, in a
coherent way with respect to the $\eta$'s and the $\mu$'s.
\goodbreak

\lem
Let $(K,\otimes,\unit)$ be a tensor triangulated category with unit. Then the
endomorphism ring $\End_K(\unit)$ is commutative.
\endlem

\proof
Let $f,g\in\End_K(\unit)$. Observe that the following diagram is commutative in
$K\times K$\,:
$$\xymatrix@C=30pt @R=30pt{(\unit,\unit)\ar[r]^-{\disp(f,\Id)}\ar[d]_-{\disp(\Id,g)}\ar[rd]|-{\disp(f,g)}
    &(\unit,\unit)\ar[d]^-{\disp(\Id,g)}
    \\
    (\unit,\unit)\ar[r]_-{\disp(f,\Id)}
    &(\unit,\unit)
}$$ Apply the functor $\otimes:K\times K\to K$ and replace $\unit\otimes\unit$
by $\unit$ using the fact that both isomorphisms $\unit\otimes\unit\to\unit$
coincide. The result follows.
\qed

\espaceavanttruc
Let $(K,\otimes,\unit)$ be a tensor triangulated category with unit. One can
define the spectrum $\Spc(K,\otimes)$ and the presheaf of triangulated
categories ${\Cal K}$ exactly as before, observing that the latter is equipped
with a unit on each open by localizing the unit of $K$. On this topological
space $\Spc(K)$, one can then replace the presheaf of rings $\End({\Cal
K})/\PNil({\Cal K})$ which we used above by the presheaf of rings $\End_{{\Cal
K}}(\unit)$. Exactly the same kind of arguments as above give us the following
result.

\thm
Consider the functor $X\mapsto(\Dperf(X),\Lotimes,{\Cal O}_X)$ from
topologically noetherian schemes to tensor triangulated categories with unit.
Then this functor is faithful and reflects isomorphisms.
\endthm
\qed

\rem
A more precise statement is true. Let ${\Cal T}$ be the category of molecular
tensor triangulated categories with unit, and geometric morphisms, as defined
in Definitions~\tref{defgeom} and~\tref{defmolecular}. Then there exists a
functor from ${\Cal T}$ into ringed spaces such that its pre-composition with
the functor $X\mapsto (\Dperf(X),\Lotimes,{\Cal O}_X)$ just gives back the
natural inclusion from topologically noetherian schemes into ringed spaces.

\goodbreak
\espaceavantpar\noindent{\bf References}
\espaceaprespar

\bible P.~Balmer, M.~Schlichting; Idempotent completion of triangulated
categories; J.\ Algebra {\bf 236}, no.~2 (2001), pp.~819-834;

\bible P.~Berthelot, A.~Grothendieck, L.~Illusie; SGA 6, Th\'eorie des
intersections et th\'eor\`eme de Riemann-Roch;Lecture Notes in Mathematics,
vol.~225, Springer 1971;

\bible A.~Bondal, D.~Orlov;Reconstruction of a variety
from the derived category and groups of autoequivalences;Compositio Math. {\bf
125}, no.~3 (2001), pp.~327-344;

\bible A.~Grothendieck, J.~A.~Dieudonn\'e; El\'ements de g\'eom\'etrie
alg\'ebrique I; Grundlehren der mathematischen Wissenschaften, Band 166,
Springer 1971;

\bible Shigeru Mukai;Duality between $D(X)$ and $D(\hat X)$ with its application to
Picard sheaves;Nagoya Math.\ J.\ {\bf 81} (1981), pp.~153-175;

\bible Amnon Neeman;The chromatic tower for D(R); with an appendix by Marcel
B\"okstedt, Topology {\bf 31}, no.~3 (1992), pp.~519-532;

\bible
Amnon Neeman;$K$-Theory for triangulated categories $3{3\over4}$: a direct
proof of the theorem of the heart;K-theory {\bf 22}, nos.~1-2 (2001),
pp.~1-144;

\bible Dmitri Orlov;Equivalences of derived categories and $K3$ surfaces;
J.\ Math.\ Sci.\ {\bf 84}, no.~5 (1997), pp.~1361--1381;

\bible Alexander Polishchuk;Symplectic biextensions and a generalization of the
Fourier-Mukai transform; Math.\ Res.\ Lett.\ {\bf 3}, no.~6 (1996),
pp.~813--828;

\bible
Jeremy Rickard;Derived categories and stable equivalence;J.\ Pure Appl.\
Algebra~{\bf 61}, no.~3 (1989), pp.~303-317;

\bible Robert W.\ Thomason;The classification of triangulated
subcategories;Compositio Math. {\bf 105}, no.~1 (1997), pp.~1-27;

\bible R.\,W.\ Thomason, T.\ Trobaugh;Higher algebraic K-theory of
schemes and of derived categories; in\,: The Grothendieck Festschrift,
vol.~III, pp.~247-435, Progr.\ Math.~88, Birk\"auser 1990;

\bible Avishay Vaknin;Virtual triangles;K-Theory {\bf 22}, nos.\ 1-2 (2001),
pp.~161-198;

\bible Jean-Louis Verdier; Des cat\'egories d\'eriv\'ees des cat\'egories ab\'eliennes;
Ast\'e\-risque {\bf 239}, So\-ci\'e\-t\'e math\'ematique de France (1996);

\bible Charles Weibel;An introduction to homological algebra;Cambridge
University Press 1994;

\vfill \footnote""{Version of April 30, 2002}

\enddocument